\documentclass[11pt]{amsart}\usepackage{latexsym,amsfonts,amsmath,amssymb,amscd}\usepackage[cm]{fullpage}
\newtheorem{Theorem}{Theorem}[section]
\newtheorem{Remark}{Remark}[section]
\newtheorem*{Question}{Question}
\newcommand{\Pf}{\noindent {\bf Proof }}
\begin{document}
\title[positive determinant]{Einstein four-manifolds with self-dual Weyl curvature of nonnegative determinant}
\author{Peng Wu}
\address{Shanghai Center for Mathematical Sciences, Fudan University, 220 Handan Road, Shanghai 200433, China} \email{wupenguin@fudan.edu.cn}\thanks{}
\dedicatory{}\date{\today}
\begin{abstract}
We prove that simply connected Einstein four-manifolds of positive scalar curvature are conformally K\"ahler if and only if the determinant of the self-dual Weyl curvature is positive.
\end{abstract}
\maketitle

\section{Introduction}

This is a sequel to the author's thesis \cite{Wuthesis} (see also \cite{Wu1}) and \cite{Wu,Wu2,Wu3}. The question that when a four-manifold with a complex structure admits a compatible Einstein metric of positive scalar curvature has been answered by Tian \cite{Tian} (see also Odaka, Spotti, and Sun \cite{OSS}), LeBrun \cite{LeBrun}, respectively. K\"ahler-Einstein four-manifolds of positive scalar curvature \cite{Tian,OSS} are $\mathbb{C}P^2$, $\mathbb{C}P^1\times\mathbb{C}P^1$, or $\mathbb{C}P^2\#k\overline{\mathbb{C}P^2}$ ($3\leq k\leq 8$). Hermitian, Einstein four-manifolds of positive scalar curvature \cite{LeBrun} are either K\"ahler-Einstein, or $\mathbb{C}P^2\#\overline{\mathbb{C}P^2}$ with Page metric \cite{Page}, or $\mathbb{C}P^2\#2\overline{\mathbb{C}P^2}$ with Chen-LeBrun-Weber metric \cite{CLW}. Recall that a Hermitian, Einstein metric is an Einstein metric which is Hermitian with respect to some integrable complex structure.

It is natural to ask, conversely,
\begin{Question}
When does a four-manifold with an Einstein metric of positive scalar curvature admit a compatible complex structure?
\end{Question}

There have been several answers to this question. A classical result of Derdzi\'nski (Theorem 2 in \cite{Der}) states that, passing to a double cover of the manifold if necessary, if the self-dual Weyl curvature $W^+$ is parallel and $\#{\rm spec}(W^+)=2$, then the metric is K\"ahler; if $\#{\rm spec}(W^+)=2$, then the metric is Hermitian, where $\#{\rm spec}(W^+)$ is the number of distinct eigenvalues of $W^+$.

Richard and Seshadri \cite{RS}, Fine, Krasnov, and Panov \cite{FKP}, and the author \cite{Wu2} proved that if the metric has half nonnegative isotropic curvature, then it is either half conformally flat or K\"ahler. LeBrun \cite{LeBrun2} proved that if $W^+(\omega,\omega)>0$ for some $\omega\in H^2_+(M)$, then the metric is Hermitian. The author \cite{Wu2} proved that if the metric has conformally half nonnegative isotropic curvature, then it is either half conformally flat or Hermitian.

The eigenvalues of $W^+$ of any K\"ahler metric on four-manifolds are $-\frac{R}{12}$, $-\frac{R}{12}$, $\frac{R}{6}$, where $R$ is the scalar curvature. LeBrun \cite{LeBrun3} proved that any Hermitian, Einstein metric of positive scalar curvature on four-manifolds must be conformal to an extremal K\"ahler metric, so the eigenvalues of $W^+$ are $-\lambda$, $-\lambda$, $2\lambda$ for some positive function $\lambda$, hence $\det W^+>0$. In this paper we prove

\begin{Theorem} \label{Thm1.1}
Simply connected Einstein four-manifolds of positive scalar curvature are conformally K\"ahler if and only if $\det W^+>0$.
\end{Theorem}

On Riemannian four-manifolds, $W^+$ is traceless, so $W^+$ satisfies a simple algebraic inequality $3\sqrt{6}|\det W^+|\leq |W^+|^3$, and the equality holds if and only if $\#{\rm spec}(W^+)\leq 2$. The idea of proving Theorem \ref{Thm1.1} is to prove that if $\det W^+>0$ then $3\sqrt{6}\det W^+\equiv |W^+|^3$, then apply the aforementioned results of Derdzi\'nski \cite{Der} and LeBrun \cite{LeBrun3}.

\begin{Remark}
According to Theorem 2 in \cite{Der}, the ``simply connected" condition in Theorem \ref{Thm1.1} can be replaced by ``oriented and $H^1(M,\mathbb{Z}_2)=0$".
\end{Remark}

The idea of the proof is motivated by previous work of Gursky and LeBrun\cite{GL}, Yang \cite{Yang}, and the author \cite{Wu} on the rigidity of Einstein four-manifolds of positive sectional curvature, in which the authors analyzed $|W^{\pm}|^2$, and reduced the problem to $W^{\pm}\equiv 0$, then applied a classical result of Hitchin (Theorem 13.30 in \cite{Besse}). As the author observed in Section 5 of \cite{Wu}, these methods might be in some sense constrained by the refined Kato inequality of Gursky and LeBrun \cite{GL}. {\it The new idea in this paper is to analyze both $|W^+|^2$ and $\det W^+$, and, instead of reducing to $W^{\pm}\equiv 0$, we reduce the problem to $3\sqrt{6}\det W^+\equiv|W^+|^3$, as explained above.}

The key step in the proof is to construct a subharmonic function of the form $f(|W^+|^2,\det W^+)$, which is based on Derdzi\'nski's derivation \cite{Der} of the Weitzenb\"ock formula for the self-dual Weyl curvature, and the author's work \cite{Wu} on an alternative proof of the refined Kato inequality, and the classification of Einstein four-manifolds of three-nonnegative curvature operator. Precisely we have

\begin{Theorem} \label{Thm1.2}
Let $(M,g)$ be a compact oriented four-manifold with $\delta W^+=0$. If $\det W^+>0$, then there exists a constant $k_0$ depending on $\min_M |W^+|^{-3}\det W^+$, $\min_M |W^+|^{-2}\det W^+$, and $\min_M R$, such that for any $k\geq k_0$,
$$F_k=|W^+|^{\frac{1}{3}}\left[1-54\bigg(\frac{\det W^+}{|W^+|^3}\bigg)^2\right]^k$$
is a subharmonic function on $M$. Furthermore by the Stokes Theorem we get that $3\sqrt{6}\det W^+\equiv|W^+|^3$.
\end{Theorem}

Interestingly, $F_k$ is closely related to the refined Kato inequality, see Remark \ref{remarkkato} in Section 2 for details.

By similar arguments we have,

\begin{Theorem} \label{Thm1.3}
Simply connected Einstein four-manifolds of positive scalar curvature and $\det W^+\geq0$ are either anti-self-dual or conformally K\"ahler.
\end{Theorem}

\begin{Theorem} \label{Thm1.4}
Compact oriented Ricci-flat four-manifolds with $H^1(M,\mathbb{Z}_2)=0$ and $\det W^+\geq0$ are anti-self-dual, therefore the universal cover of $M$ is either $\mathbb{R}^4$ with flat metric or a K3 surface with Calabi-Yau metric.
\end{Theorem}

Theorem 1.1 and its proof suggest us to ask the following question,
\begin{Question}
Are simply connected Einstein four-manifolds of positive scalar curvature conformally K\"ahler, if the self-dual Weyl curvature is nonvanishing?
\end{Question}

We would like to point out that recently LeBrun \cite{LeBrun4} gave an alternative proof of Theorem \ref{Thm1.1} based on the method in \cite{LeBrun2}. Furthermore he relaxed the condition in Theorem \ref{Thm1.1} to $W^+\neq 0$ and $|W^+|^{-3}\det W^+\geq-\frac{5\sqrt{2}}{21\sqrt{21}}$.

\begin{Remark}
We observe that on Einstein four-manifolds of positive scalar curvature, either $W^{\pm}\equiv 0$ or the average of $\det W^{\pm}$ has a positive lower bound. Recall the Weitzenb\"ock formula of Derdzi\'nski \cite{Besse,Der},
\begin{equation*}
\begin{split}
\Delta|W^{\pm}|^2=2|\nabla W^{\pm}|^2+R|W^{\pm}|^2-36\det W^{\pm}.
\end{split}
\end{equation*}
In our paper, we use $\Delta f={\rm tr}\nabla^2 f=g^{ij}\nabla_i\nabla_j f$ for $f\in C^{\infty}(M)$. Gursky and LeBrun \cite{GL} proved that either $W^{\pm}\equiv 0$ or $\int_M|W^{\pm}|^2dv\geq\int_M\frac{R^2}{24}dv$. Combining the two formulas together we get, either $W^{\pm}\equiv 0$, or
\begin{equation*}
\begin{split}
\int_M \det W^{\pm}dv\geq 2\int_M\frac{R^3}{12^3}dv.
\end{split}
\end{equation*}
\end{Remark}

\textbf{Acknowledgement.} The author thanks his advisors Professors Xianzhe Dai and Guofang Wei for their guidance, encouragement, and constant support.  The author thanks Professors Claude LeBrun and Yuan Yuan for helpful discussions. The author thanks the anonymous referee for many suggestions that greatly improve the presentation of the paper. The author acknowledges the hospitality of Capital Normal University, East China Normal University, and Mathematical Sciences Research Institute, where part of this work was carried out. The author was partially supported by NSFC No.11701093 and China recruit program for global young talents.

\section{Proof}

We explain the method of constructing subharmonic functions of the form $f(|W^+|^2,\det W^+)$ on $M$ in two steps.

Step 1. We briefly recall Derdzi\'nski's derivation of the Weitzenb\"ock formula for Riemannian metrics of $\delta W^{\pm}=0$ on four-manifolds.

Let $\lambda_1\leq\lambda_2\leq\lambda_3$ be the eigenvalues of $W^+$, with corresponding orthogonal eigenvectors
$$\omega_1=e^1\wedge e^2+e^3\wedge e^4,\quad \omega_2=e^1\wedge e^3+e^4\wedge e^2,\quad \omega_3=e^1\wedge e^4+e^2\wedge e^3,$$
then $W^+$ can be expressed as
$$W^+=\frac{1}{2}(\lambda_1 \omega_1\otimes \omega_1+\lambda_2 \omega_2\otimes \omega_2+\lambda_3 \omega_3\otimes \omega_3).$$

Let $M_W$ be the open dense subset of $M$, consisting of points at which the number of distinct eigenvalues of $W^+$ is locally constant, then $\lambda_i$ and $\omega_i$ ($i=1,2,3$) may be assumed differentiable in a neighborhood of any point $p\in M_W$, so there exist 1-forms $a,b,c$ defined near $p$, such that
\begin{equation*}
\begin{split}
\nabla \omega_1=&a\otimes \omega_2-c\otimes \omega_3,\\
\nabla \omega_2=&b\otimes \omega_3-a\otimes \omega_1,\\
\nabla \omega_3=&c\otimes \omega_1-b\otimes \omega_2.
\end{split}
\end{equation*}

By analyzing the Ricci identities for $\omega_1$, $\omega_2$, $\omega_3$, Derdzi\'nski proved that, if $\delta W^+=0$, then in a neighborhood of $p\in M_W$,
\begin{equation*}
\begin{split}
\nabla\lambda_1=&(\lambda_2-\lambda_1)(\iota_{a^{\#}}\omega_3)^{\#}+(\lambda_3-\lambda_1)(\iota_{c^{\#}}\omega_2)^{\#},\\
\nabla\lambda_2=&(\lambda_1-\lambda_2)(\iota_{a^{\#}}\omega_3)^{\#}+(\lambda_3-\lambda_2)(\iota_{b^{\#}}\omega_1)^{\#},\\
\nabla\lambda_3=&(\lambda_1-\lambda_3)(\iota_{c^{\#}}\omega_2)^{\#}+(\lambda_2-\lambda_3)(\iota_{b^{\#}}\omega_1)^{\#},\\
\Delta \lambda_1=&2(\lambda_1-\lambda_2)|(\iota_{a^{\#}}\omega_3)^{\#}|^2+2(\lambda_1-\lambda_3)|(\iota_{c^{\#}}\omega_2)^{\#}|^2
+\frac{R}{2}\lambda_1-2\lambda_1^2-4\lambda_2\lambda_3,\\
\Delta \lambda_2=&2(\lambda_2-\lambda_1)|(\iota_{a^{\#}}\omega_3)^{\#}|^2+2(\lambda_2-\lambda_3)|(\iota_{b^{\#}}\omega_1)^{\#}|^2
+\frac{R}{2}\lambda_2-2\lambda_2^2-4\lambda_1\lambda_3,\\
\Delta \lambda_3=&2(\lambda_3-\lambda_1)|(\iota_{c^{\#}}\omega_2)^{\#}|^2+2(\lambda_3-\lambda_2)|(\iota_{b^{\#}}\omega_1)^{\#}|^2
+\frac{R}{2}\lambda_3-2\lambda_3^2-4\lambda_1\lambda_2.
\end{split}
\end{equation*}
where $\iota$ is the interior product, $\#$ is the sharp operator.

\begin{Remark}
Derdzi\'nski also derived the formula for $\nabla W^+$, combining these formulas together he proved the classical Weitzenb\"ock formula,
\begin{equation*}
\begin{split}
\Delta |W^+|^2=2|\nabla W^+|^2+R|W^+|^2-36\det W^+.
\end{split}
\end{equation*}
\end{Remark}

Step 2. We reduce the subharmonicity of functions of the form $f(|W^+|^2,\det W^+)$ on $M$ to a system of partial differential inequalities on $\mathbb{R}^2$, based on the author's alternative proof of the refined Kato inequality.

There are only two nontrivial elementary symmetric polynomials of $\lambda_1,\lambda_2,\lambda_3$:
\begin{equation*}
\begin{split}
\sigma_1=&\lambda_1+\lambda_2+\lambda_3=0,\\
x\triangleq-2\sigma_2=&\lambda_1^2+\lambda_2^2+\lambda_3^2=|W^+|^2,\\
y\triangleq\sigma_3=&\lambda_1\lambda_2\lambda_3=\det W^+.
\end{split}
\end{equation*}

For simplicity, we define vector fields $X\triangleq(\lambda_1-\lambda_2)(\iota_{a^{\#}}\omega_3)^{\#}$, $Y\triangleq(\lambda_2-\lambda_3)(\iota_{b^{\#}}\omega_1)^{\#}$, $Z\triangleq(\lambda_3-\lambda_1)(\iota_{c^{\#}}\omega_2)^{\#}$ in a neighborhood of $p\in M_W$. We have
\begin{equation*}
\begin{split}
\nabla x=&-2(\lambda_1-\lambda_2)X-2(\lambda_2-\lambda_3)Y-2(\lambda_3-\lambda_1)Z,\\
\nabla y=&\lambda_3(\lambda_1-\lambda_2)X+\lambda_1(\lambda_2-\lambda_3)Y+\lambda_2(\lambda_3-\lambda_1)Z,\\
\Delta x=&8|X|^2+8|Y|^2+8|Z|^2-4\langle X,Y\rangle-4\langle X,Z\rangle-4\langle Y,Z\rangle+(Rx-36y),\\
\Delta y=&-4\lambda_3|X|^2-4\lambda_1|Y|^2-4\lambda_2|Z|^2-4\lambda_2\langle X,Y\rangle-4\lambda_1\langle X,Z\rangle-4\lambda_3\langle Y,Z\rangle\\ +&\Big(\frac{3}{2}Ry-x^2\Big).
\end{split}
\end{equation*}

Let $f=f(x,y)$ be a differentiable function on $M$. On $M_W$, we have
\begin{equation}\label{Deltaf}
\begin{split}
\Delta f=&f_x\Delta x+f_y\Delta y+f_{xx}|\nabla x|^2+f_{yy}|\nabla y|^2+2f_{xy}\nabla x\nabla y\\
=&[8f_x-4\lambda_3f_y+(\lambda_1-\lambda_2)^2(4f_{xx}+\lambda_3^2f_{yy}-4\lambda_3f_{xy})]|X|^2\\
+&[8f_x-4\lambda_1f_y+(\lambda_2-\lambda_3)^2(4f_{xx}+\lambda_1^2f_{yy}-4\lambda_1f_{xy})]|Y|^2\\
+&[8f_x-4\lambda_2f_y+(\lambda_3-\lambda_1)^2(4f_{xx}+\lambda_2^2f_{yy}-4\lambda_2f_{xy})]|Z|^2\\
+&2[-2f_x-2\lambda_2f_y+(\lambda_1-\lambda_2)(\lambda_2-\lambda_3)(4f_{xx}+\lambda_1\lambda_3f_{yy}+2\lambda_2f_{xy})]\langle X,Y\rangle\\
+&2[-2f_x-2\lambda_1f_y+(\lambda_1-\lambda_2)(\lambda_3-\lambda_1)(4f_{xx}+\lambda_2\lambda_3f_{yy}+2\lambda_1f_{xy})]\langle X,Z\rangle\\
+&2[-2f_x-2\lambda_3f_y+(\lambda_2-\lambda_3)(\lambda_3-\lambda_1)(4f_{xx}+\lambda_1\lambda_2f_{yy}+2\lambda_3f_{xy})]\langle Y,Z\rangle\\
+&(Rx-36y)f_x+\Big(\frac{3}{2}Ry-x^2\Big)f_y.
\end{split}
\end{equation}

We denote $A$, $B$, $C$ as the coefficients of $|X|^2$, $|Y|^2$, $|Z|^2$ in Equation (\ref{Deltaf}), respectively; and $2D$, $2E$, $2F$ as the coefficients of $\langle X,Y\rangle$, $\langle X,Z\rangle$, $\langle Y,Z\rangle$ in Equation (\ref{Deltaf}), respectively. We define
\begin{equation*}
\begin{split}
{\rm\bf I}\triangleq& A|X|^2+B|Y|^2+C|Z|^2+2D\langle X,Y\rangle+2E\langle X,Z\rangle+2F\langle Y,Z\rangle,\\
{\rm\bf II}\triangleq& (Rx-36y)f_x+\Big(\frac{3}{2}Ry-x^2\Big)f_y.
\end{split}
\end{equation*}

Then we have
\begin{equation*}
\begin{split}
\Delta f=&{\rm\bf I}+{\rm\bf II}.
\end{split}
\end{equation*}

If ${\rm\bf I}\geq0$ and ${\rm\bf II}\geq 0$ on $M_W$, then $\Delta f\geq 0$ on $M_W$, moreover since $M_W$ is an open dense subset of $M$ and $f$ is differentiable, we conclude that $\Delta f\geq0$ on $M$.

\

We consider ${\rm\bf I}$ as a quadratic form of (components of) $X,Y,Z$. In order for ${\rm\bf I}\geq 0$, we need $A>0, B>0, C>0$. Consider ${\rm\bf I}$ as a quadratic function of (components of) $X$, then its minimum is
\begin{equation*}
\tilde{{\rm\bf I}}=A^{-1}[(AB-D^2)|Y|^2+(AC-E^2)|Z|^2+2(AF-DE)\langle Y,Z\rangle].
\end{equation*}

In order for $\tilde{{\rm\bf I}}\geq 0$, we need $AB-D^2>0,AC-E^2>0$. Consider $\tilde{{\rm\bf I}}$ as a quadratic function of (components of) $Y$, then its minimum is
\begin{equation*}
(AB-D^2)^{-1}(ABC-AF^2-BE^2-CD^2+2DEF)|Z|^2.
\end{equation*}

Therefore the quadratic form ${\rm\bf I}\geq 0$ if $A,B,C,D,E,F$ satisfy the following system
\begin{equation*}
\begin{cases}
\hspace{-0.3cm}&A>0,\quad B>0,\quad C>0,\\
\hspace{-0.3cm}&\mathrm{I}_{31}\triangleq AB-D^2>0,\\
\hspace{-0.3cm}&\mathrm{I}_{32}\triangleq AC-E^2>0,\\
\hspace{-0.3cm}&\mathrm{I}_{33}\triangleq BC-F^2>0,\\
\hspace{-0.3cm}&\mathrm{I}_4\triangleq ABC-AF^2-BE^2-CD^2+2DEF\geq 0.
\end{cases}
\end{equation*}

Notice that the characterization of the quadratic form ${\rm\bf I}>0$ follows from Sylvester's criterion.

Observe that $A+B+C>0$, $\mathrm{I}_{31}+\mathrm{I}_{32}+\mathrm{I}_{33}>0$, and $\mathrm{I}_4\geq 0$ will ensure that all $A$, $B$, $C$, $\mathrm{I}_{31}$, $\mathrm{I}_{32}$, $\mathrm{I}_{33}$, are positive. Therefore $f$ is a subharmonic function on $M$ if
\begin{equation*}
\begin{cases}
\hspace{-0.3cm}&0\leq{\rm I}_1\triangleq {\rm\bf II}=(Rx-36y)f_x+(\frac{3}{2}Ry-x^2)f_y;\\
\hspace{-0.3cm}&0<{\rm I}_2\triangleq A+B+C;\\
\hspace{-0.3cm}&0<{\rm I}_3\triangleq AB+AC+BC-D^2-E^2-F^2;\\
\hspace{-0.3cm}&0\leq{\rm I}_4= ABC-AF^2-BE^2-CD^2+2DEF.
\end{cases}
\end{equation*}

Plugging in $A$, $B$, $C$, $D$, $E$, $F$ to the above system, we conclude that $f(x,y)$ is a subharmonic function on $M$, if $f(x,y)$, considering as a function on $\Omega=\{(x,y)\in\mathbb{R}^2:\ x^3\geq 54y^2\}\subset\mathbb{R}^2$, satisfies the following system of partial differential inequalities,
\begin{equation*}
({\rm PDI})\begin{cases}
0\leq&\hspace{-0.3cm}\mathrm{I}_1=(Rx-36y)f_x+(\frac{3}{2}Ry-x^2)f_y.\\
0<&\hspace{-0.3cm}\mathrm{I}_2=24xf_{xx}+72yf_{xy}+x^2f_{yy}+48f_x.\\
0<&\hspace{-0.3cm}\mathrm{I}_3=6(x^3-54y^2)(f_{xx}f_{yy}-f_{xy}^2)+24(7xf_x-12yf_y)f_{xx}\\
\hspace{0.3cm}&\hspace{0.125cm}+\,8(63yf_x-2x^2f_y)f_{xy}+x(7xf_x-12yf_y)f_{yy}+180f_x^2-12xf_y^2.\\
0\leq&\hspace{-0.3cm}\mathrm{I}_4=6(x^3-54y^2)(f_{xx}f_{yy}-f_{xy}^2)f_x+4(30xf_x^2-72yf_xf_y-x^2f_y^2)f_{xx}\\
\hspace{0.3cm}&\hspace{0.125cm}+\,4(90yf_x^2-4x^2f_xf_y-3xyf_y^2)f_{xy}+(5x^2f_x^2-12xyf_xf_y-9y^2f_y^2)f_{yy}\\
\hspace{0.3cm}&\hspace{0.125cm}+\,100f_x^3-14xf_xf_y^2-8yf_y^3.
\end{cases}
\end{equation*}

\Pf of Theorem \ref{Thm1.2}. We will construct a function $f(x,y)$ that satisfies Sytem (PDI) in the subregion $\Omega_{\delta}=\{y\geq\delta>0,x^3\geq54y^2\}\subset\Omega$.

Define $z=x^{-\frac{3}{2}}y\in [-\frac{1}{3\sqrt{6}},\frac{1}{3\sqrt{6}}]$ at points where $x\neq 0$, and $f(x,y)=x^{\frac{1}{6}}(1-54z^2)^k$, plugging $f_x$, $f_y$, $f_{xx}$, $f_{xy}$, $f_{yy}$ into System (PDI), we have
\begin{equation*}
\begin{split}
\mathrm{I}_1=&\frac{1}{6}x^{\frac{1}{6}}(1-54z^2)^k\big[R+36(18k-1)x^{-1}y\big],\\
\mathrm{I}_2=&\frac{2}{3}x^{-\frac{5}{6}}(1-54z^2)^{k-1}\big[54(18k-1)(18k+7)z^2-(162k-7)\big],\\
\mathrm{I}_3=&\frac{2}{9}x^{-\frac{5}{3}}(1-54z^2)^{2k-2}\big[2916(18k-1)^2(18k+5)z^4-108(1944k^2-162k+5)z^2\\
&-(162k-5)\big],\\
\mathrm{I}_4\equiv&0.
\end{split}
\end{equation*}

Since $M$ is compact and $x^3\geq 54y^2$, if $y\geq\delta$ for some $\delta>0$, then $z\geq\delta_1$, $x^{-1}y\geq\delta_2$, for some $\delta_1>0$, $\delta_2>0$. By choosing $k$ large enough, we get that ${\rm I}_i\geq0$, moreover ${\rm I}_i>0$ when $1-54z^2>0$, $i=1,2,3$. So we have ${\rm\bf I}\geq0$ and ${\rm\bf II}\geq0$ on $M_W$, therefore $\Delta f\geq0$ on $M$.

By Stokes Theorem we get $\Delta f\equiv 0$ on $M$, then ${\rm\bf I}\equiv 0$, ${\rm\bf II}\equiv 0$ on $M_W$. From ${\rm\bf II}\equiv 0$ on $M_W$ we get that $1-54z^2\equiv 0$ on $M_W$, which implies $3\sqrt{6}y\equiv x^{\frac{3}{2}}$ on $M$,

\qed

\textbf{Proof} of Theorem \ref{Thm1.1}. By Theorem \ref{Thm1.2} and the aforementioned results in \cite{Der,LeBrun3}, $(M,g)$ is conformally K\"ahler.

\qed

\begin{Remark} \label{remarkkato}
Recall the refined Kato inequality \cite{GL} for $W^+$ of Einstein metrics on four-manifolds,
$$|\nabla W^+|^2\geq\frac{5}{3}|\nabla|W^+||^2.$$
Consider a function $f(x)=f(|W^+|^2)$, by the Weitzenb\"ock formula, we have
\begin{equation*}
\begin{split}
\Delta f=&f'\Delta x+f''|\nabla x|^2\\
=&2|\nabla W^+|^2f'+4|\nabla|W^+||^2xf''+(Rx-36y)f'.
\end{split}
\end{equation*}
Denote $\Delta_D f=2|\nabla W^+|^2f'+4x|\nabla|W^+||^2f''$, the ``derivative part" of the Weitzenb\"ock formula, then the refined Kato inequality for $W^+$ can be interpreted as
$$\Delta_D x^{\frac{1}{6}}\geq 0.$$
Moreover, $\frac{1}{6}$ is the smallest power such that this inequality holds, see Section 5 in \cite{Wu} for details. The function we construct, $x^{\frac{1}{6}}(1-54z^2)^k$, can be considered as a homogeneous variation of $x^{\frac{1}{6}}$, since $z$ depends only on the quotient $\frac{\lambda_1}{\lambda_3}$, but is independent of the magnitude of $W^+$.
\end{Remark}

\

\textbf{Proof} of Theorem \ref{Thm1.3}. We will construct a function $f(x,y)$ that satisfies System (PDI) in the subregion $\Omega_0=\{x>0,y\geq 0,x^3\geq54y^2\}\subset\Omega$.

Consider $f=x^{\frac{1}{6}}h(z)$, $z\in[0,\frac{1}{3\sqrt{6}}]$, with $h(z)\geq0$ to be determined. Plugging $f_x$, $f_y$, $f_{xx}$, $f_{xy}$, $f_{yy}$ into System (PDI), we have
\begin{equation*}
\begin{split}
\mathrm{I}_1=&\frac{1}{6}x^{-1}\big[(Rx-36y)h-6x^{\frac{3}{2}}(1-54z^2)h'\big].\\
\mathrm{I}_2=&\frac{1}{3}x^{-\frac{5}{6}}\big[3(1-54z^2)h''-270zh'+14h\big].\\
\mathrm{I}_3=&\frac{1}{9}x^{-\frac{5}{3}}\big[-27z(1-54z^2)h'h''+3(1-54z^2)hh''-6(2-351z^2)h'^2\\
&-324zhh'+10h^2\big].\\
\mathrm{I}_4\equiv& 0.
\end{split}
\end{equation*}

Suppose $h'(z)=(1-54z^2)^{-1}\phi(z)h(z)$ with $\phi(z)$ to be determined, then we have
\begin{equation*}
\begin{split}
h''=&(1-54z^2)^{-2}\big[(1-54z^2)\phi'+\phi^2+108z\phi\big]h.
\end{split}
\end{equation*}

Plugging into the above system, we have
\begin{equation*}
\begin{split}
\mathrm{I}_1=&\frac{1}{6}hx^{\frac{1}{6}}[R-6x^{\frac{1}{2}}(\phi+6z)].\\
\mathrm{I}_2=&\frac{1}{3}hx^{-\frac{5}{6}}(1-54z^2)^{-1}\big[3(1-54z^2)\phi'+3\phi^2+54z\phi+14(1-54z^2)\big].\\
\mathrm{I}_3=&\frac{1}{9}h^2x^{-\frac{5}{3}}(1-54z^2)^{-2}\big[3(1-54z^2)(1-54z^2-9z\phi)\phi'-27z\phi^3\\
&-9(1+108z^2)\phi^2+10(1-54z^2)^2\big].
\end{split}
\end{equation*}

First notice that $\max_{\Omega_0}x^{\frac{1}{2}}=\infty$, so ${\rm I}_1\geq0$ in $\Omega_0$ if and only if $\phi+6z\leq0$. It is obvious that if $\phi+6z\leq0$, then $1-54z^2-9z\phi\geq0$ when $z\geq0$.

Next notice that $(1-54z^2-9z\phi)h\mathrm{I}_2-\mathrm{I}_3=12\phi^2h^2+4(1-54z^2-9z\phi)^2h^2$, so if ${\rm I}_3\geq0$ in $\Omega_0$ and $\phi+6z\leq0$ then ${\rm I_2}\geq0$ in $\Omega_0$.

In summary, if ${\rm I}_3\geq0$ and $\phi+6z\leq0$, then ${\rm I}_1\geq0$ and ${\rm I}_2\geq0$. ${\rm I}_3\geq0$ and $\phi+6z\leq0$ is equivalent to an Abel differential inequality of the second kind \cite{Murphy} on $[0,\frac{1}{3\sqrt{6}}]$ with a constraint condition,
\begin{equation}\label{abelforphi}
\begin{split}
3(1-54z^2)(1-54z^2-9z\phi)\phi'-27z\phi^3-9(1+108z^2)\phi^2&\\
+10(1-54z^2)^2&\geq 0,\\
\phi+6z&\leq 0.
\end{split}
\end{equation}

To further simplify the system, we denote $\psi(z)=\phi(z)+6z$, then the constrained Abel differential inequality (\ref{abelforphi}) can be written as
\begin{equation} \label{abelforpsi}
\begin{split}
3(1-54z^2)(1-9z\psi)\psi'-27z\psi^3-9(1+54z^2)\psi^2&\\
+270z\psi-8(1+54z^2)&\geq 0,\\
\psi&\leq 0.
\end{split}
\end{equation}

We choose the initial value $\psi(0)=-6\sqrt{6}$, then the constrained Abel differential inequality (\ref{abelforpsi}) has a solution $\psi(z)$ on $[0,\frac{1}{3\sqrt{6}}]$, which is monotonically increasing and $\psi(\frac{1}{3\sqrt{6}})=-3\sqrt{6}$.  By the definition of $\psi(z)$, one can check that $h(z)$ is monotonically decreasing on $[0,\frac{1}{3\sqrt{6}}]$, and $h(\frac{1}{3\sqrt{6}})=0$. So we have ${\rm I}_i\geq 0$, moreover ${\rm I}_i>0$ when $1-54z^2>0$, $i=1,2,3$. Therefore ${\rm I}_1\geq0$, ${\rm I}_2\geq0$ on $M_W\backslash\{W^+=0\}$, and $\Delta f\geq0$ on $M_W\backslash\{W^+=0\}$.

Furthermore, by the definition of $\phi$, the function $f$ we construct satisfies the property that $f^{6k}=x^kh^{6k}(z)\in C^2(M)$ for sufficient large $k$. By the above argument, we have $\Delta f^{6k}\geq0$ on $M$. By Stokes Theorem, we get $\Delta f^{6k}\equiv 0$, then ${\rm\bf I}\equiv 0$, ${\rm\bf II}\equiv 0$ on $M_W$. From ${\rm\bf II}=0$ on $M_W$, we get that either $x=0$ or $1-54z^2=0$ on $M_W$, therefore, either $x=0$, or $x^3=54y^2$ and $y>0$ on $M$. By Prop 5 in \cite{Der}, we get that either $x\equiv0$, or $x^3\equiv54y^2$ and $y>0$ on $M$. Therefore $(M,g)$ is either anti-self-dual or conformally K\"ahler.

\qed

\textbf{Proof} of Theorem \ref{Thm1.4}. Notice that the nonnegativity of ${\rm\bf I}$ is independent of the sign of the scalar curvature. By the same argument as in the proof of Theorem \ref{Thm1.3}, if $R=0$, then ${\rm I}_1=-hx^{\frac{2}{3}}(\phi+6z)=-hx^{\frac{2}{3}}\psi$. Therefore $\psi$ satisfies the same Abel differential inequality with the same constraint condition, and we get the same conclusion that either $x\equiv0$, or $x^3\equiv 54y^2$ and $y>0$ on $M$.

If $x^3\equiv54y^2$ and $y>0$, then $(M,g)$ is conformally K\"ahler. By Proposition 5 in \cite{Der}, $\bar g=(24x)^{\frac{1}{3}}g$ is a K\"ahler metric with scalar curvature $\bar R=(24x)^{\frac{1}{6}}>0$. On the other hand, by the conformal change of the scalar curvature, $\bar R$ has to be nonpositive somewhere, which leads to a contradiction.

Therefore we have $x\equiv 0$, that is, $(M,g)$ is anti-self-dual.

\qed

\begin{Remark} \label{remarkPDE}
It is interesting to observe in the proof of Theorem \ref{Thm1.3} that all homogeneous variations of $x^{\frac{1}{6}}$, that is, functions of the form $x^{\frac{1}{6}}h(z)$, where $h(z)$ is an arbitrary differentiable function, solve the partial differential equation
\begin{equation*}
\begin{split}
0={\rm I}_4&=6(x^3-54y^2)(f_{xx}f_{yy}-f_{xy}^2)f_x+4(30xf_x^2-72yf_xf_y-x^2f_y^2)f_{xx}\\
&+4(90yf_x^2-4x^2f_xf_y-3xyf_y^2)f_{xy}+(5x^2f_x^2-12xyf_xf_y-9y^2f_y^2)f_{yy}\\
&+100f_x^3-14xf_xf_y^2-8yf_y^3.
\end{split}
\end{equation*}

One may ask whether this equation admits solutions of a different form, which may help us to characterize K\"ahler-Einstein or Hermitian, Einstein metrics using different curvature conditions by constructing functions $f(x,y)$ that satisfies System (PDI) in different subsets of $\Omega$.
\end{Remark}

\end{document}